\documentclass[12pt,twoside]{article}
\usepackage{amssymb}
\font\teneufm=eufm10 scaled \magstep1
\font\seveneufm=eufm7 scaled \magstep1
\font\fiveeufm=eufm5  scaled \magstep1
\newfam\eufmfam
\textfont\eufmfam=\teneufm
\scriptfont\eufmfam=\seveneufm
\scriptscriptfont\eufmfam=\fiveeufm

\usepackage{mathrsfs}

\newfam\msbfam
\font\tenmsb=msbm10 scaled \magstep1  \textfont\msbfam=\tenmsb
\font\sevenmsb=msbm7 scaled \magstep1 \scriptfont\msbfam=\sevenmsb
\font\fivemsb=msbm5 scaled \magstep1  \scriptscriptfont\msbfam=\fivemsb
\def\Bbb{\fam\msbfam \tenmsb}

\def\RR{{\Bbb R}}
\def\CC{{\Bbb C}}

\def\ZZ{{\Bbb Z}}
\def\TT{{\Bbb T}}

\def\ra{\rightarrow}

 \def\HollowBoxx #1#2#3{{\dimen0=#1 \advance\dimen0 by -#2
       \dimen1=#1 \advance\dimen1 by #3
        \vrule height 0pt depth #3 width #2
       \hskip -#3
       \vrule height #1 depth #3 width #3}}
 \def\LeftContraction{\mathord{\kern1.45pt \HollowBoxx{6pt}{3.5pt}{.4pt}}\,}

 \def\HollowBox #1#2#3{{\dimen0=#1 \advance\dimen0 by -#3
       \dimen1=#1 \advance\dimen1 by #3
        \vrule height #1 depth #3 width #3
        \vrule height 0pt depth #3 width #2
        \hskip -#3}}
 \def\RightContraction{\mathord{\, \HollowBox{6pt}{3.1pt}{.4pt}} \kern1.6pt}

\def\qed{{\hfill $\Box$}}
\newtheorem{theorem}{THEOREM}[section]
\newtheorem{corollary}[theorem]{Corollary}

\newtheorem{remark}[theorem]{Remark}

\begin{document}

\begin{center}
{\Large \bf A Remark on a Theorem
\medskip\\
by Kodama and Shimizu}
\medskip \\
\normalsize A. V. Isaev
\end{center}

\begin{quotation} \small \sl We prove a characterization theorem for the unit polydisc $\Delta^n\subset\CC^n$ in the spirit of a recent result due to Kodama and Shimizu. We show that if $M$ is a connected $n$-dimensional complex manifold such that (i) the group $\hbox{Aut}(M)$ of holomorphic automorphisms of $M$ acts on $M$ with compact isotropy subgroups, and (ii) $\hbox{Aut}(M)$ and $\hbox{Aut}(\Delta^n)$ are isomorphic as topological groups equipped with the compact-open topology, then $M$ is holomorphically equivalent to $\Delta^n$.
\end{quotation}

\thispagestyle{empty}

\pagestyle{myheadings}
\markboth{A. V. Isaev}{ A Remark on a Theorem by Kodama and Shimizu}

\setcounter{section}{0}

\section{Introduction}
\setcounter{equation}{0}

For a connected complex manifold $M$, let $\hbox{Aut}(M)$ denote the group of holomorphic automorphisms of $M$. Endowed with the compact-open topology, $\hbox{Aut}(M)$ is a topological group. We are interested in characterizing complex manifolds by their automorphism groups. 

In general, two complex manifolds $M_1$ and $M_2$ need not be holomorphically equivalent if the topological groups $\hbox{Aut}(M_1)$ and $\hbox{Aut}(M_2)$ are isomorphic. A simple example of this kind with non-trivial automorphism groups is given by spherical shells
$$
S_r:=\{z\in\CC^n: r<||z||<1\}, \quad 0\le r<1.
$$
It is straightforward to see that for $n\ge 2$ the group $\hbox{Aut}(S_r)$ coincides with the unitary group $U_n$ for all $r$. Next, every $S_r$ is a Kobayashi-hyperbolic Reinhardt domain. It is shown in \cite{Kr}, \cite{S} that two such domains are holomorphically equivalent if and only if they are equivalent by means of an elementary algebraic map, i.e. a map of the form
$$
z_j\mapsto\lambda_j z_1^{a_{j1}}\cdot\dots\cdot z_n^{a_{jn}},\quad j=1,\dots,n,
$$
where $\lambda_j\in\CC^*$ and $a_{jk}$ are integers satisfying $\hbox{det}(a_{jk})\ne 0$. An elementary algebraic map is holomorphic and one-to-one on $S_r$ only if it is linear (i.e. reduces to dilations and a permutation of coordinates). However, $S_{r_1}$ and $S_{r_2}$ are not equivalent by means of such a linear map for $r_1\ne r_2$.

If the group $\hbox{Aut}(M)$ is sufficiently large, one can hope to obtain positive characterization results. For example, it was shown in \cite{IK} that the space $\CC^n$ is completely characterized by its holomorphic automorphism group as follows: if $M$ is a connected complex manifold of dimension $n$ and the groups $\hbox{Aut}(M)$ and $\hbox{Aut}(\CC^n)$ are isomorphic as topological groups, then $M$ is holomorphically equivalent to $\CC^n$. A similar characterization was obtained for the unit ball $B^n\subset\CC^n$ in \cite{I} (see also the erratum) and, under certain additional assumptions (that will be discussed below), for direct products $B^k\times\CC^{n-k}$ in \cite{BKS} as well as for the space $\CC^n$ without some coordinate hyperplanes in \cite{KS1}, \cite{KS2}.

Recently, in \cite{KS3} Kodama and Shimizu obtained the following characterization of another classical domain, the unit polydisc $\Delta^n\subset\CC^n$ (the direct product of $n$ copies of the unit disc $\Delta\subset\CC$).

\begin{theorem}\label{ks3}\rm[KS3]\,\sl Let $M$ be a connected complex manifold of dimension $n$ that is holomorphically separable and admits a smooth envelope of holomorphy. If $\hbox{Aut}(M)$ and $\hbox{Aut}(\Delta^n)$ are isomorphic as topological groups, then $M$ is holomorphically equivalent to $\Delta^n$.
\end{theorem}
In particular, Theorem \ref{ks3} holds for Stein manifolds and for all domains in $\CC^n$.

The connected component of the identity $\hbox{Aut}(\Delta^n)^0$ of the group $\hbox{Aut}(\Delta^n)$ is isomorphic to the direct product of $n$ copies of the group $\hbox{Aut}(\Delta)\simeq SU_{1,1}/\ZZ_2$, and therefore contains a subgroup (which is a maximal compact subgroup) isomorphic to the $n$-torus $\TT^n$. A topological group isomorphism between $\hbox{Aut}(M)$ and $\hbox{Aut}(\Delta^n)$ yields a smooth action by holomorphic transformations of $\TT^n$ on $M$. The assumptions of holomorphic separability and smoothness of the envelope of holomorphy in Theorem \ref{ks3} are used by the authors to linearize this action thus representing the manifold $M$ as a Reinhardt domain in $\CC^n$. This is possible due to a theorem by Barrett, Bedford and Dadok (see \cite{BBD}). We note that similar assumptions were imposed on manifolds in \cite{BKS}, \cite{KS1}, \cite{KS2} to guarantee the applicability of the result of \cite{BBD}.

It is anticipated that the assertion of Theorem \ref{ks3} remains true if the assumptions of holomorphic separability and smoothness of the envelope of holomorphy are dropped. In this note we offer a version of Theorem \ref{ks3} in this direction. In particular, we do not refer to the linearization result of \cite{BBD} in our proofs. Instead, we require that for every $p\in M$ the isotropy subgroup
$$
\hbox{Aut}_p(M):=\left\{g\in\hbox{Aut}(M): g(p)=p\right\}
$$
is compact in $\hbox{Aut}(M)$ and linearize the action of $\hbox{Aut}_p(M)$ near $p$, which is possible due to the results of Bochner in \cite{B} (see also \cite{Ka}). We note that the linearizability of actions of compact groups on complex manifolds with fixed points goes back to H. Cartan (see \cite{M} for an account of Cartan's results of this kind). In fact, we will only use the faithfulness of the isotropy representation (defined below); this statement is known as Cartan's uniqueness theorem (see \cite{C}). The local linearizability (as opposed to the global linearizability of the $\TT^n$-action) is sufficient to characterize $\Delta^n$. It is not clear at this time how one could avoid using linearization arguments altogether. One difficulty here is the low-dimensionality of the maximal compact subgroup of $\hbox{Aut}(\Delta^n)^0$. For comparison, the maximal compact subgroup of $\hbox{Aut}(B^n)$ is isomorphic to $U_n$ and thus has dimension $n^2$. This fact was of great help in \cite{I} (see also \cite{IK}).  

Our result is the following theorem.

\begin{theorem}\label{main}\sl Let $M$ be a connected complex manifold of dimension $n$ such that for every $p\in M$ the isotropy subgroup $\hbox{Aut}_p(M)$ is compact in $\hbox{Aut}(M)$. If $\hbox{Aut}(M)$ and $\hbox{Aut}(\Delta^n)$ are isomorphic as topological groups, then $M$ is holomorphically equivalent to $\Delta^n$.
\end{theorem}

We remark that the assumption of compactness of the isotropy subgroups holds for large classes of manifolds a priori not covered by Theorem \ref{ks3}. For example, it holds whenever the action of the group $\hbox{Aut}(M)$ on the manifold $M$ is proper, i.e. the map
$$
\hbox{Aut}(M)\times M\ra M\times M,\quad (g,p)\mapsto (g(p),p)
$$
is proper. It is shown in \cite{Ka} that $\hbox{Aut}(M)$ acts on $M$ properly if and only if one can find a continuous $\hbox{Aut}(M)$-invariant distance on $M$. In particular, the action of $\hbox{Aut}(M)$ is proper for all Kobayashi-hyperbolic manifolds (see also \cite{Ko}). Hence the following holds (cf. Remark \ref{rem}).

\begin{corollary}\label{cor}\sl Let $M$ be a connected Kobayashi-hyperbolic manifold of dimension $n$. If $\hbox{Aut}(M)$ and $\hbox{Aut}(\Delta^n)$ are isomorphic as topological groups, then $M$ is holomorphically equivalent to $\Delta^n$.
\end{corollary}

\section{Proof of Theorem \ref{main}}

Let $\hbox{Aut}(M)^0$ be the connected component of the identity of $\hbox{Aut}(M)$. Since $\hbox{Aut}(\Delta^n)^0$ is a Lie group of dimension $3n$ in the compact-open topology, so is $\hbox{Aut}(M)^0$. Furthermore, every maximal compact subgroup of $\hbox{Aut}(M)^0$ is $n$-dimensional and isomorphic to $\TT^n$. For every $p\in M$ the subgroup $\hbox{Aut}_p(M)^c:=\hbox{Aut}_p(M)\cap\, \hbox{Aut}(M)^0$ is compact and therefore is contained in some maximal compact subgroup of $\hbox{Aut}(M)^0$. Since the dimension of the $\hbox{Aut}(M)^0$-orbit of $p$ cannot exceed $2n$, it follows that $\hbox{dim}\,\hbox{Aut}_p(M)^c=n$. Hence $\hbox{Aut}_p(M)^c$ is a maximal compact subgroup of $\hbox{Aut}(M)^0$ (thus $\hbox{Aut}_p(M)^c=\hbox{Aut}_p(M)^0$), and the action of $\hbox{Aut}(M)^0$ on $M$ is transitive.       

Let
$$
\alpha_p:\, \hbox{Aut}_p(M)^0\ra GL(\RR,T_p(M)),\quad g\mapsto dg(p)
$$
be the isotropy representation of $\hbox{Aut}_p(M)^0$, where $T_p(M)$ is the tangent space to $M$ at $p$ and $dg(p)$ is the differential of a map $g$ at $p$. Let further
$$
L_p:=\alpha_p\left(\hbox{Aut}_p(M)^0\right)
$$
be the corresponding linear isotropy subgroup. By the results of \cite{C}, \cite{B}, \cite{Ka} the isotropy representation is continuous and faithful. In particular, $L_p$ is a compact subgroup of $GL(\RR,T_p(M))$ isomorphic to $\hbox{Aut}_p(M)^0$. In some coordinates in $T_p(M)$ the group $L_p$ becomes a subgroup of the unitary group $U_n$. Since $L_p$ is isomorphic to $\TT^n$, it is conjugate in $U_n$ to the subgroup of all diagonal unitary matrices. In particular, for every $p\in M$ the group $L_p$ contains the element $-\hbox{id}$.

Let ${\mathscr G}$ be an $\hbox{Aut}(M)^0$-invariant Hermitian metric on $M$. Since $\hbox{Aut}(M)^0$ acts on $M$ transitively, such a metric can be constructed by choosing an $L_{p_0}$-invariant positive-definite Hermitian form on $T_{p_0}(M)$ for some $p_0\in M$, and by extending it to a Hermitian metric on all of $M$ using the $\hbox{Aut}(M)^0$-action (see \cite{P} for the existence of invariant metrics for not necessarily transitive proper actions). The manifold $M$ equipped with the metric ${\mathscr G}$ is a Hermitian symmetric space.

The theorem now follows from the general theory of Hermitian symmetric spaces (see \cite{H}). Indeed, since the group $\hbox{Aut}(M)^0$ acts on $M$ with compact isotropy subgroups, contains a symmetry at every point of $M$, is semi-simple and is isomorphic to the direct product of $n$ copies of the simple group $SU_{1,1}/\ZZ_2$, the manifold $M$ is holomorphically isometric to the product of $n$ one-dimensional irreducible Hermitian symmetric spaces (see Theorem 3.3 in Chapter IV, Theorems 1.1 and 4.1 in Chapter V, Propositions 4.4, 5.5 and Theorem 6.1 in Chapter VIII of \cite{H}). Clearly, each of the one-dimensional irreducible Hermitian symmetric spaces must be equivalent to the unit disc $\Delta$, and the proof is complete. \qed   

\begin{remark}\label{rem}\rm One can obtain Corollary \ref{cor} without  referring to the theory of Hermitian symmetric spaces. Indeed, as in the proof of Theorem \ref{main}, we see that $M$ is homogeneous. Hence, by the (non-trivial) result of \cite{N}, the manifold $M$ is holomorphically equivalent to a bounded domain in $\CC^n$. Corollary \ref{cor} now follows from Theorem \ref{ks3}.
\end{remark}

{\obeylines
Department of Mathematics
The Australian National University
Canberra, ACT 0200
AUSTRALIA
E-mail: alexander.isaev@maths.anu.edu.au
}

\end{document}